\documentclass[letterpaper, 10 pt, conference]{ieeeconf}  % Comment this line out
\IEEEoverridecommandlockouts                              % This command is only
\usepackage{graphicx}
\usepackage[utf8]{inputenc}
\usepackage[T1]{fontenc}
\usepackage{algpseudocode}
\usepackage{algorithm}
\usepackage{multicol}
\usepackage{tikz}
\usepackage{siunitx}
\usepackage{amsmath,amsfonts,amssymb}
\usepackage{mathtools}
\usetikzlibrary{automata,arrows,positioning,calc}
\usepackage{romannum}
\usepackage{tikz}
\usepackage{float}
\usepackage{url}

\usepackage{algorithm, algpseudocode}
\usepackage{lipsum}


\newtheorem{definition}{Definition}

\newtheorem{remark}{Remark}

\title{\LARGE \bf
A Log-domain Interior Point Method for Convex Quadratic Games
}

\author{Bingqi Liu$^1$ and Dominic Liao-McPherson$^1$% <-this % stops a space
\thanks{$^1$The University of British Columbia, Vancouver, BC, V6T 1Z4 Canada. Emails: \texttt{bingqil@students.ubc.ca}, \texttt{dliaomcp@mech.ubc.ca}.}
\thanks{This work was supported by the Natural Sciences and Engineering Research Council of Canada (Reference Number RGPIN-2023-03257).}%
}

\begin{document}

\maketitle
\thispagestyle{empty}
\pagestyle{empty}

%%%%%%%%%%%%%%%%%%%%%%%%%%%%%%%%%%%%%%%%%%%%%%%%%%%%%%%%%%%%%%%%%%%%%%%%%%%%%%%%
\begin{abstract}
In this paper, we propose an equilibrium-seeking algorithm for finding generalized Nash equilibria of non-cooperative monotone convex quadratic games. Specifically, we recast the Nash equilibrium-seeking problem as variational inequality problem that we solve using a log-domain interior point method and provide a general purpose solver based on this algorithm. This approach is suitable for non-potential, general sum games and does not require extensive structural assumptions. We demonstrate the efficiency and versatility of our method using two benchmark games and demonstrate our algorithm is especially effective on small to medium scale problems.
\end{abstract}

%%%%%%%%%%%%%%%%%%%%%%%%%%%%%%%%%%%%%%%%%%%%%%%%%%%%%%%%%%%%%%%%%%%%%%%%%%%%%%%%
\section{INTRODUCTION}
Many important problems in engineering involve collections of interacting self-interested agents, including network routing, power grid management, market analysis, supply chain optimization, etc. \cite{c1,c2,c3,c4,c5}. Non-cooperative games are a powerful tool for modelling these complex systems; in particular, the Nash equilibria of Convex Quadratic Games (CQGs) can model important applications in various resource-constrained scenarios, e.g., electric vehicle charging, traffic network routing, and social welfare distribution \cite{c1,c3,c9}.

Roughly speaking, a CQG is a collection of coupled quadratic programming (QP) problems. There are a variety of methods for finding Generalized (i.e., constrained) Nash Equilibria (GNE's) of CQGs; these can be roughly classified into iterative best response (BR) methods and VI-type methods. BR methods are intuitive, they seek an equilbirum through repreated play \cite{c8} but are not guaranteed to converge \cite{c10}, especially for problems with coupling constraints.

VI-type methods are the basis for many successful algorithms \cite{c1,c9,c11,c12}. Those methods rely on the fact that useful subset of GNEs must satisfy a variational inequality \cite{c8,c13} which can be solved using a variety of methods. The state-of-the-art for CQGs are built on first-order operator-theoretic methods e.g., the forward-backward splitting (FBS) or the Douglas-Rashford splitting (DRS), and result in provably convergent algorithms that are well suited to distributed settings and parallelization \cite{c1,c9,c14}. However, these approaches are “bespoke” in that they generally require extensive reformulation of the problem and make strong assumptions about its properties, e.g., (strict) monotonicity, potential, or aggregative structure. 

Another approach is to use second-order methods (e.g., interior point sequential quadratic programming or augmented Lagrangian methods \cite{c8,c11,c18}) to solve the VI. Second-order methods are generally more difficult to distribute and/or parallelize but typically converge faster and require fewer assumptions about the problem structure. They can solve a much broader class of problems without requiring extensive reformulation including non-monotone problems. However, existing implementations \cite{c12,c15} are specialized to nonlinear problems with a specific multi-stage structure (essentially multi-agent robotic trajectory planning problems).

In this paper, we present an algorithm for finding variational GNEs (v-GNEs) of CQGs based on a recently proposed log-domain interior point method (IPM) for convex quadratic programming \cite{c17} and demonstrate the utility of our method through comparisons with OS methods on three benchmark problems. Our goal is to provide a general purpose method (like those available for QPs) that can be used to solve a wide-variety of problems with little need for extensive reformulations (as can be the case when using operator splitting methods). 

The paper is structured as follows: In Section II, we derive a method for finding GNEs of a CQG using a log-domain IPM, then give a detailed description of the algorithm that we proposed to solve QGs. Section III provided three classical examples of how our method can be applied. Section IV evaluates the performance of our method on different problem difficulties. In Section V, we conclude and suggest future research directions.

\subsection{Notation}
The set of first n natural numbers $\{1,2,3,...,n\}$ is denoted by $[n]$. The perpendicular sign $\perp$ denotes the complementary: defining $\odot$ as the element-wise product, $\nu \geq 0 \perp s\geq0$ implies $\nu \odot s = 0$, $s \geq 0$, $\nu \geq 0$. We denote the vertical concatenation of N matrices or vectors by $(A_1,A_2,..,A_N)$. Without special notation, we assume vectors are column vectors. Given N vectors $x_1,...,x_N$, $\mathbf{x}:=(x_1,...,x_N)$, $\mathbf{\Bar{x}}:=\frac{1}{N}\sum_{i=1}^{N}x_i$, and $x_{-i}:=(x_1,...,x_{i-1},x_{i+1},...,x_N)$. $\mathbf{0}(\mathbf{1})$ denotes a vector with all elements equal to $0\,(1)$. $I_n$ and $0_n$ denote the square matrix of size $n\times n$ with diagonal elements equal to 1 and 0 respectively. $I$ denotes an identity matrix with suitable size. For a given vector $a\in\mathbb{R}^n$, $e^a\in\mathbb{R}^n$ and $\forall i \in n$, ${e^a}_i=e^{a_i}$.

\subsection{Source code}
The source code of the IPM CQG solver is available at: \url{https://github.com/BL-DLM-paper/CVX-Game-solver}.

\section{Convex Quadratic Games}
Consider a group of N players indexed by $i\in[N]$ with decision variable $x_i \in \mathbb{R}^n$. The aim of the $i^{th}$ player is to solve the following QP:
\begin{subequations} \label{eq:CQG}
\begin{align}
\min_{x_i} &~ \frac{1}{2}x_i^\intercal Q_i x_i+p_i^\intercal x_i + \sum_{i\neq j}x_i^\intercal S_{ij}x_j\\
\textrm{s.t.} &~~ A\mathbf{x}+b \geq \mathbf{0}\\
&~~ G\mathbf{x}+h=\mathbf{0},
\end{align}
\end{subequations}
where $Q_i\succeq 0$, $S_{ij}$, $p_i$ are problem data and the set
\begin{equation}
    \chi = \{\mathbf{x} \in \mathbb{R}^{nN}|A\mathbf{x}+b \geq \mathbf{0},~~ G\mathbf{x}+h = \mathbf{0} \}
\end{equation}
encodes the constraints. These constraints can be local or global depending on the structure of the problem data $A$, $G$, $b$ and $h$. 

A strategy $\mathbf{\tilde{x}}=(\tilde{x}_1,...,\tilde{x}_N)$ is a GNE of the game (1) if
\begin{equation}
J_i(\tilde{x}_i,\tilde{x}_{-i})\leq \min_{x_i}\{J_i(x_i,\tilde{x}_{-i})~|~(x_i,\tilde{x}_{-i}) \in \chi\}
\end{equation}
where $J_i(x,x_i) =\frac{1}{2}x_i^\intercal Q_i x_i+p_i^\intercal x_i + \sum_{i\neq j}x_i^\intercal S_{ij}x_j$. Since $Q_i \succeq 0$, each $J_i$ is convex in the local decision variable $x_i$ and the VI
\begin{equation} \label{eq:VI}
(x-\tilde x)^T (W\tilde x + f) \geq 0,\, \forall x \in \chi
\end{equation}
where
\begin{equation} \label{eq:W-def}
W=\begin{bmatrix}
Q_1 &S_{12} &S_{13} &\hdots &S_{1N}\\
S_{21} &Q_2 &S_{23} &\hdots &S_{2N}\\
& & \ddots & &\\
& & & \ddots &\\
S_{N1} &S_{N2} &S_{N3} &\hdots &Q_N\\
\end{bmatrix},
f=\begin{bmatrix}
p_1\\
p_2\\
\vdots\\ 
\vdots\\
p_N\\
\end{bmatrix}.
\end{equation}
is necessary for $\tilde x$ to be a GNE of the CQG \eqref{eq:CQG} \cite{c19}. If $\chi \neq \emptyset$ so that Slater's condition holds the VI \eqref{eq:VI} is equivalent to the KKT conditions
\begin{subequations} \label{eq:KKT}
\begin{align}
W\mathbf{x}+f-A^\intercal \nu -G^\intercal \lambda = \mathbf{0} \\
\mathbf{0}\leq \nu \perp A\mathbf{x}+b \geq \mathbf{0} \label{eq:comp}\\
G\mathbf{x}+h=\mathbf{0}.
\end{align}
\end{subequations}
GNEs that satisfy \eqref{eq:VI} are known as variational GNE's (v-GNEs). In addition to computational tractability, they have a variety of desirable properties including a notion of ``larger social stability'', see \cite[Section 5]{c20}. However, they are a subset of all GNEs of \eqref{eq:CQG}; in particular a game might have several GNEs but no v-GNEs. Luckily, the existence of a v-GNE can be guaranteed under some conditions.

\begin{definition}
The monotonicity constant $\sigma$ of the CQG \eqref{eq:CQG} is $\sigma = \lambda_{min}(\frac{W+W^\intercal}{2})$ where $\lambda_{min}$ denotes the smallest eigenvalue.
\end{definition}

If $\sigma \geq 0$ and $\chi \neq \emptyset$ then the KKT conditions \eqref{eq:KKT} are necessary and sufficient for a v-GNE of \eqref{eq:CQG} \cite[Theorem 1]{c9} and a v-GNE is guaranteed to exist \cite[Prop 12.11]{c21}.

\begin{remark}
If the Hessian matrix W of the CQG is symmetric then \eqref{eq:KKT} are the KKT conditions of a quadratic program and finding a v-GNE of \eqref{eq:CQG} reduces to an optimization problem (i.e., \eqref{eq:CQG} is a potential game). If $W \neq W^\intercal$ then standard optimization algorithms cannot be applied.
\end{remark}

\section{A Log-domain Interior Point Method}
In this section, we outline a log-domain IPM to solve \eqref{eq:KKT} for a v-GNE of the CQG. IPMs are homotopy methods that smooth the complimentarity condition \eqref{eq:comp} with a variable smoothing parameter $\mu > 0$ such that the modified KKT system is easy to solve for large $\mu$ and we recover a solution of the original problem as $\mu \to 0$. IPMs approximately trace a homotopy path (known as the ``central path'') towards a solution as $\mu \to 0$.

\subsection{The Central Path}
The log-domain IPM \cite{c17} is a recently proposed variation of IPM that uses exponential functions to smooth the complementarity conditions. We first reformulate \eqref{eq:KKT} by adding the slack variable $s = Ax + b$, this leads to the equations
\begin{equation}
\begin{aligned}
W\mathbf{x}+f-A^\intercal \nu -G^\intercal \lambda = \mathbf{0}\\
\mathbf{0}\leq \nu \perp s\geq \mathbf{0}\\
s = Ax + b\\
G\mathbf{x}+h=\mathbf{0}.
\end{aligned}
\end{equation}
Then we define the smoothed equations by introducing the auxilliary variable $v$ and barrier parameter $\mu > 0$ and making the substitutions $\nu = \sqrt{\mu} e^{-v}$ and $s = \sqrt{\mu} e^v$, this leads to the ``log-domain central path'' equations
\begin{equation} \label{eq:central-path}
\begin{aligned}
W\mathbf{x}+f-A^\intercal \sqrt{\mu}e^{v} -G^\intercal \lambda = \mathbf{0}\\
A\mathbf{x}+b = \sqrt{\mu}e^{-v}\\
G\mathbf{x}+h=\mathbf{0}.
\end{aligned}
\end{equation}
As $\mu \to 0$, solutions of \eqref{eq:central-path} approach the solution of the original VI \eqref{eq:VI} \cite[Lemme 3.3]{c17} and thus a v-GNE of the original game. The main advantage of this log-domain homotopy over typical interior-point approaches is that this parameterization enforces the conditions $s\geq 0$, $\nu \geq 0$ by construction.

\subsection{Damped Newton's Method}
In practice, the log-domain IPM approximately tracks solutions of the central-path equations \eqref{eq:central-path} for decreasing $\mu$ using a damped Newton’s method. The core of the algorithm is the Newton iteration
\begin{equation}
    v^+ = v + \alpha \Delta v
\end{equation}
where $\alpha$ is a step size chosen to enforce convergence, $\Delta v$ is computed by solving the Newton-step system
\begin{multline} \label{eq:newton-step}
    \begin{bmatrix}
        W  & - \sqrt{\mu} A^{\intercal} U(v)  & -G^\intercal\\
        A  &  \sqrt[]{\mu} U(-v)  & 0\\
        G  & 0  & 0\\
    \end{bmatrix}
    \begin{bmatrix}
        x^+\\
        \Delta v\\
        \lambda^+\\
    \end{bmatrix} \\
    =\begin{bmatrix}
        \sqrt{\mu}A^{\intercal} e^v-f\\
        \sqrt{\mu}e^{-v}-b\\
        -h\\
    \end{bmatrix}
\end{multline}
and $U(v)=diag\{e^v\}$. The algorithm alternates between centering steps that attempt to solve the central path equations \eqref{eq:central-path} for a fixed $\mu$ and reduction steps that decrease $\mu$. The step length is selected so as to ensure convergence of the underlying Newton's method and $\mu$ is reduced so as to remain within the quadratic convergence region of Newton's method (more details are available in \cite{c17}). The method is described formally in Algorithm~1 and is globally convergent under the condition that $\lambda_{min}(\frac{W+W^\intercal}{2}+A^\intercal A) > 0$ \cite{c17}.

\begin{algorithm} \label{algo:IPM}
\caption{Log-domain IPM}
\begin{algorithmic}[1]
\Procedure{IPM}{$W$, $A$, $G$, $f$, $b$, $h$, $\varepsilon$, $\tau$}
    \State Set $\mu=1$, and $v=1$
    \While{$||\Delta v||_{\infty}>1 \hspace{1mm} or \hspace{1mm} \mu > \varepsilon $}
    \State $\mu^+$ = LINE-SEARCH($\mu$, $v$)
    \State [$v^+$, $x^+$, $\lambda^+$] = NEWTON($\mu^+$, $v$, $\eta_1$)
    \State $v$, $x$, $\lambda$, $\mu$ = $v^+$, $x^+$, $\lambda^+$, $\mu^+$
    \EndWhile  \label{roy's loop}
    \State [$v^+$, $x^+$, $\lambda^+$] = NEWTON($\mu$, $v$, $\eta_2$)
    \State return $x^+$
\EndProcedure

\Procedure{newton}{$\mu$, $v$, $\eta$}
    \State $i=1$
    \While{$i \leq \eta$}
    \State [$x^+$, $\Delta v$, $\lambda^+$, $\alpha$] = NEWTON-STEP($\mu$, $v$)
    \State $v^+$ = $v$ + $\alpha \Delta v$
    \State $i = i + 1$
    \If{$||\Delta v|| \leq \tau$}
        \State $v=v^+$
        \State break
    \EndIf
    \State $v=v^+$
    \EndWhile
    \State $v^+ = v$
    \State return ($x^+$, $v^+$, $\lambda^+$)
\EndProcedure

\Procedure{line-search}{$\mu$, $v$}
    \State $\mu^+$ = $\min$($\mu$, $\inf$\{$\mu > 0$ : $||\Delta v(v,\mu)||_{\infty}\leq 1$\})
    \State return $\mu^+$
\EndProcedure

\Procedure{newton-direction}{$\mu$, $v$}
    \State Solve (7) for $x^+$, $\Delta v$, $\lambda^+$
    \State $\alpha$ = $1/\max$\{$1$, $\frac{1}{2\beta} {||\Delta v||}^2_\infty$\}
    \State return ($x^+$, $\Delta v$, $\lambda^+$, $\alpha$)
\EndProcedure
\end{algorithmic}
\end{algorithm}

\subsection{Implementation Details}
Our algorithm takes in the matrices $W,A, G$ and vectors $f, b, h$ as problem data in the form of sparse matrices along with two user defined precision levels $\varepsilon , \tau> 0$. To improve numerical performance we rescale the problem data. For $W$ and $f$ we set $W \gets W/(||f||_{\infty}+1)$ and $f\gets f/(||f||_{\infty}+1)$. For $A,b,G$ and $h$ we scale each row e.g., $A_i=\frac{A_i}{|b_i|+1}$ and $b_i=\frac{b_i}{|b_i|+1}$ where $A_i$ and $b_i$ are the $i^{th}$ row respectively. The same process for $G$ and $h$ is analogous.

When solving the Newton-step system \eqref{eq:newton-step} we eliminate the second row through block elimination and introduce a diagonal regularization term $\gamma \approx 10^{-9}$ to improve numerical conditioning, this results in the following system 
\begin{multline}
    \begin{bmatrix}
        W+A^\intercal U(2v) A & -G^\intercal\\
        G  & \gamma I\\
    \end{bmatrix}
    \begin{bmatrix}
        x^+\\
        \lambda^+\\
    \end{bmatrix}\\
    =\begin{bmatrix}
        2\sqrt{\mu}A^{\intercal} e^v-f-A^\intercal U(2v)b\\
        -h\\
    \end{bmatrix}
\end{multline}
that we solve with a Sparse LU factorization routine. The factorization can be reused whenever $v$ isn't changed, this happens several times in the algorithm.

\section{Case Studies} \label{ss:case-studies}
In this section we describe three test problems from the literature that we will use to benchmark our solver. We will use the following residual of the VI to quantify distance from the central path:
\begin{equation}
R(\mathbf{x},v,\mu) =\left\|
\begin{bmatrix}
            W\mathbf{x}+f-\sqrt{\mu}A^\intercal e^v \\
            A\mathbf{x}+b+\sqrt{\mu}A^\intercal e^{-v} \\
            G\mathbf{x}+h\\
        \end{bmatrix}
\right\|.
\end{equation}

\subsection{Example 1: An electric vehicle charging game}
Consider an EV fleet with N vehicles indexed by $i\in[N]$ as in \cite{c9}. The state-of-charge of vehicle $i$ at time $t\in [n]$ is denoted by $q_i(t)$ and evolves according to $q_i(t+1) = q_i(t) + \eta_i x_i(t)$ where $\eta > 0$ is the charging efficiency and $x_i(t)$ is the control input.

The $i^{th}$ vehicle buys energy at the price
\begin{equation}
C_i \Bar{\mathbf{x}}(t) + C_i c(t)
\end{equation}
where the $C_i$ is the individual price factor for $i^{th}$ vehicle, $\bar{\mathbf{x}}(t) = 1/N\sum_{i\in [N]} x_i(t)$ is the average charging load and $c(t)$ is the non-EV energy demand at time $t$. The cost function for each vehicle is then
\begin{multline}
    J_i(x_i,x_{-i}) = \sum_{t\in [n]} \frac{1}{2} x_i(t)^\intercal Q_i x_i(t) + p_i^\intercal x_i \\+{(C_i \Bar{\mathbf{x}} (t) + C_i c)}^\intercal x_i(t)
\end{multline}
where the $Q_i$ and $p_i$ terms penalize battery degradation. 

The EVs are subject to the local constraints $\forall t\in [n]$
\begin{gather}
    0\leq x_i(t)\leq \tilde{x}_i(t),\\
     \sum_{t=1}^{n} x_i(t) \geq l_i
\end{gather}
that enforce the charging limit $\tilde{x}_i(t) > 0$ and ensure that the vehicles reaches a minimum state of charge $l_i$. Simultanously the coupling constraint
\begin{equation}
\frac{1}{N} \sum_{i=1}^{N} x_{i,t}\leq K_t,\, \forall t\in[n]
\end{equation}
where the $K(t) > 0$ is the total charging power cap which enforces a limit on the total power available from the grid.

This results in a CQG of the form \eqref{eq:CQG} which is monotone under some mild conditions on the cost functions (e.g., $Q_i \succeq 0$ and $C = (C_i)_{i\in [N]}$ is monotone \cite[Lemma 6]{c9}). The parameters $\{C_1, C_2,...,C_N\}$ determine whether the problem is symmetric (and thus a potential game). We use an sequence $C_i=1+ki$, with $k=0.01$, to ensure the game is non-potential and generate all other parameters randomly as in \cite[Section V]{c9} for a scenario where there is a spike in the background demand as depicted in Fig~\ref{fig:fig1}. The detailed derivation of the KKT system can be found in the Appendix.

\begin{figure}[thpb]
  \centering
  \includegraphics[scale=0.15]{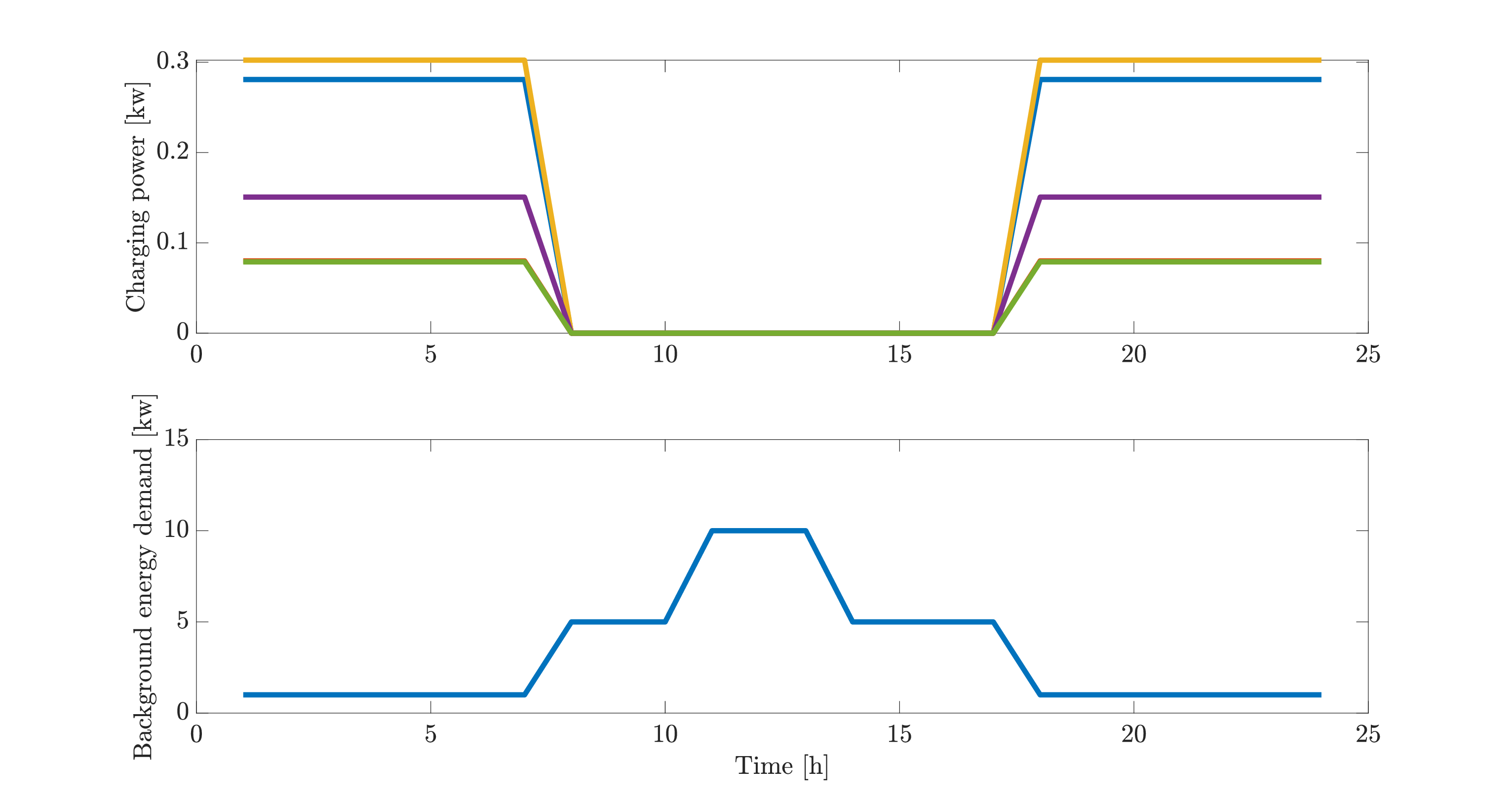}
  \caption{Non-EV power demand and EV charging power for a sample EV charging game.}
  \label{fig:fig1}
\end{figure}

\begin{figure}[thpb]
      \centering
      \includegraphics[scale=0.15]{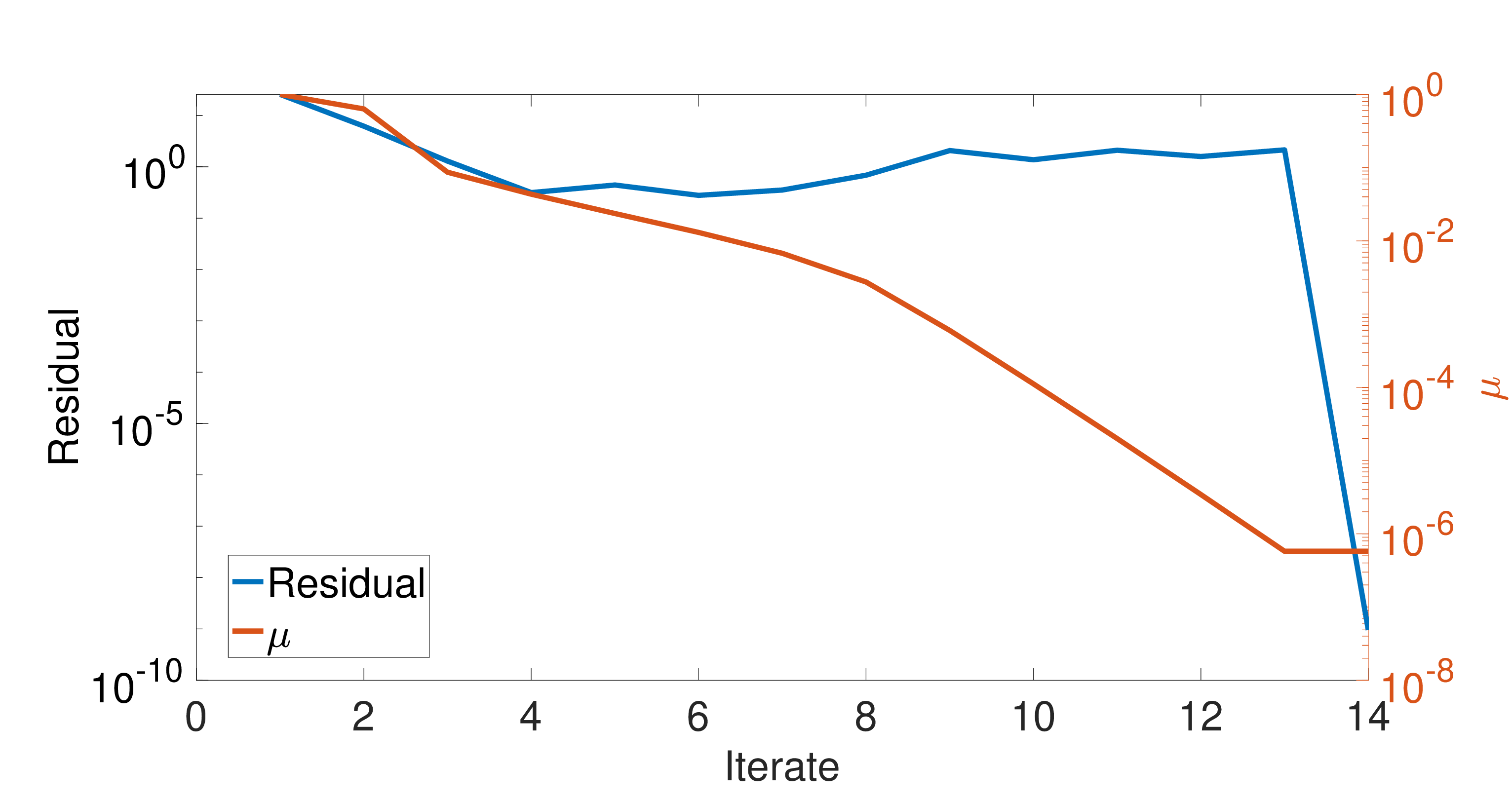}
      \caption{VI residual and barrier parameter $\mu$ at each iteration when solving the sample game. The orange curve denotes the value of $\mu$, and the blue curve denotes the value of residual. It demonstrates that the $\mu$ remains constant during the additional Newton steps while the residual is greatly reduced}
      \label{fig:fig2}
   \end{figure}

Fig~\ref{fig:fig1} shows a representative solution of the changing game. As expected, the EVs all avoid charging during the peak power demand period
The convergence behaviour of the algorithm typically follows a two-phase process which can be seen in Fig~\ref{fig:fig2}. In the first phase the algorithm tracks the central path loosely, focusing on reducing the barrier parameter $\mu$ at the expense of decreasing the KKT residual. Once $\mu$ is sufficiently small, the algorithm rapidly reduces the residual using Newton’s method to obtain a precise solution.

% if $\lambda_{min}(\frac{W+W^\intercal}{2})\geq0$ this game satisfies $A^\intercal A + W$ is monotone. We provide the details of how the concatenated KKT system suitable for our solver is derived in Section A of the Appendix. All the problem parameters have been generated randomly as Section V Part B of [9].

 % at $t$ charging interval. (10) forms the coupling constraints for this game.

% The local strategy set $\Omega_i$ is
% % \begin{equation}
% \begin{align}
% \Omega_i:=\{x_i\in \mathbb{R}^n|\ \, \forall t\in[n] \}\\
% \end{align}
% \end{equation}
% which makes sure that at every time step the $i^{th}$ vehicle's charging power won't exceed its individual upper charging input bound $\tilde{x}_i(t)$ and that at the end of $T$, it reaches a desired state of charge $l_i$.

% The final local objective function for the $i^{th}$ EV is
% \begin{equation}
% J_i(x_i,\Bar{\mathbf{x}})=\frac{1}{2} x_i^\intercal Q_i x_i+p_i^\intercal x_i +{(C_i \Bar{\mathbf{x}}+c_i)}^\intercal x_i.\\
% \end{equation}

%  We define $c=(d(1),...,d(n))$ in which $d(t)$ is the non-EV energy consumption rate at $t^{th}$ time step. The greater the total charging power demand from the EV becomes, the higher the energy price will be.

\subsection{Example 2: A market game}
Competition between companies producing similar products can be modelled with a with a market game where companies optimize the allocation of their production capacity and sales among locations to maximize profit \cite{c1}. Consider $i\in [N]$ companies (players) producing the same product that compete over $l\in[m]$ regions. Each company can set their production and sales levels at every location l, denoted by $g_{I,l}$ and $s_{I,l}$ respectively. The cost of production at location $l$ is
\begin{equation}
f_{i,l}(g_{i,l})=q_{i,l}g_{i,l}^2+g_{i,l}p_{i,l}
\end{equation}
with $q_{i,l}$ and $p_{i,l}$ are scaling parameters for the $i^{th}$ company. 
 
We adopt a simple market model in which the price is an affine function of market demand and average sales in the market. The sales price at location $l$ is
\begin{equation}
p_l(\mathbf{\Bar{s}_l})=d_l-C_i\mathbf{\Bar{s}_l}
\end{equation}
where $d_l>0$ is the total demand for the product at location $l$, $C_i>0$ is a price factor specific to the product produced by the $i^{th}$ company according to its reputation, and $\Bar{s}_l=\frac{1}{N} \sum_{i=1}^{N}{s_{i,l}}$.
The decision variables of $i^{th}$ company are $x_i=(g_{i,1},...,g_{i,m},s_{i,1},...,s_{i,m}) \in \mathbb{R}^{2m}$.
Each company would like to maximize profit, leading to the objective function for the $i^{th}$ company
\begin{equation}
J_i(x_i,\Bar{\mathbf{x}})=x_i^\intercal Q_i x_i+p_i^\intercal x_i + C_i(\Delta \Bar{\mathbf{x}})^\intercal x_i
\end{equation}
where $Q_i:=diag(q_{i,1},...,q_{i,m},0,...,0)$, $\Delta=diag(\mathbf{0},I_m)$, and $p_i:=(p_{i,1},...,p_{i,m},-d_1,...,-d_m)$

The local constraints are
\begin{gather}
    \sum_{l=1}^{m}{g_{i,l}}\geq \sum_{l=1}^{m}{s_{i,l}},\\
    g_{i,l},s_{i,j}\geq 0, g_{i,l} \leq u_{i,l}, \forall l\in [m]
\end{gather}
which ensures that each company's production level is greater than the sales level while the production level and sales level are all within a reasonable range. We add an additional coupling constraints that the total production at location $l$ must meet the demand $d_l$ while obeying a maximum capacity constraint that encodes environmental and energy consumption limits,
\begin{equation}
d_l\leq \sum_{i=1}^{N}{g_{i,l}}\leq r_l, \, \forall l\in[m].
\end{equation}

\begin{figure}[thpb]
      \centering
      \includegraphics[scale=0.2]{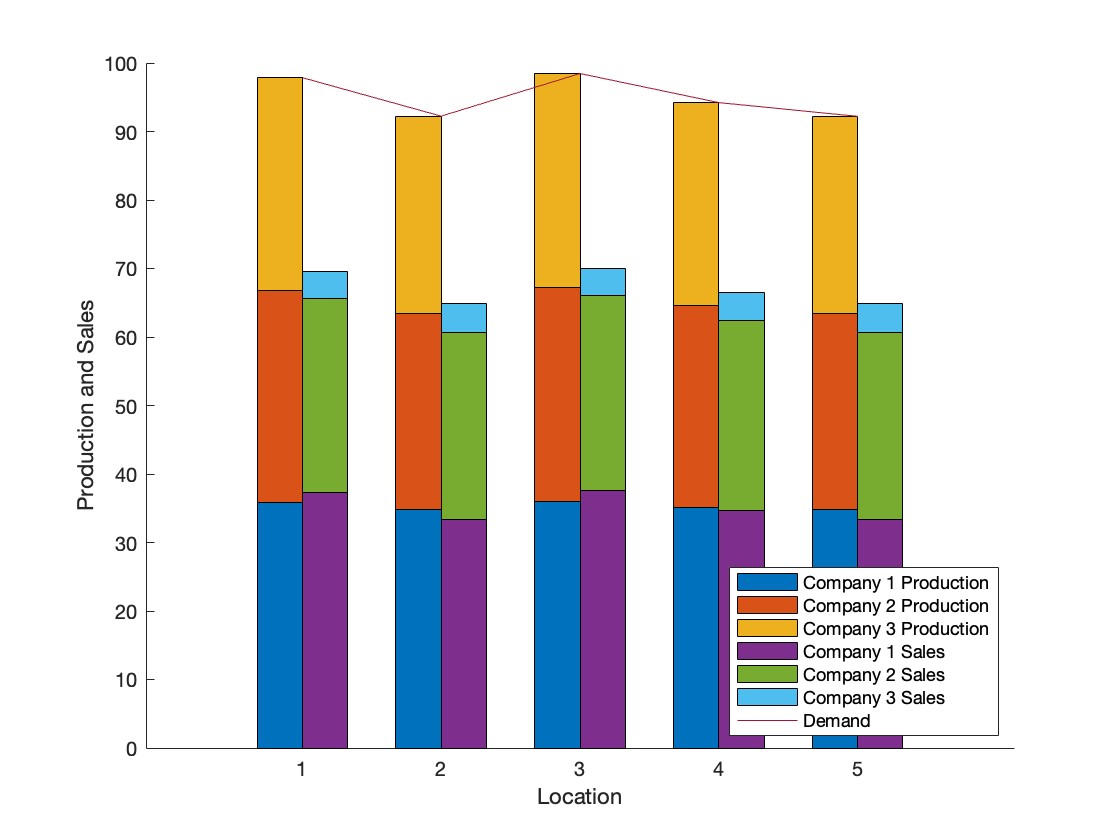}
      \caption{Demonstration of the solution to the sample market game. The production and sales level at each location is represented in stacked bar plot along with the demand level.}
      \label{fig:fig3}
   \end{figure}

\begin{figure}[thpb]
      \centering
      \includegraphics[scale=0.15]{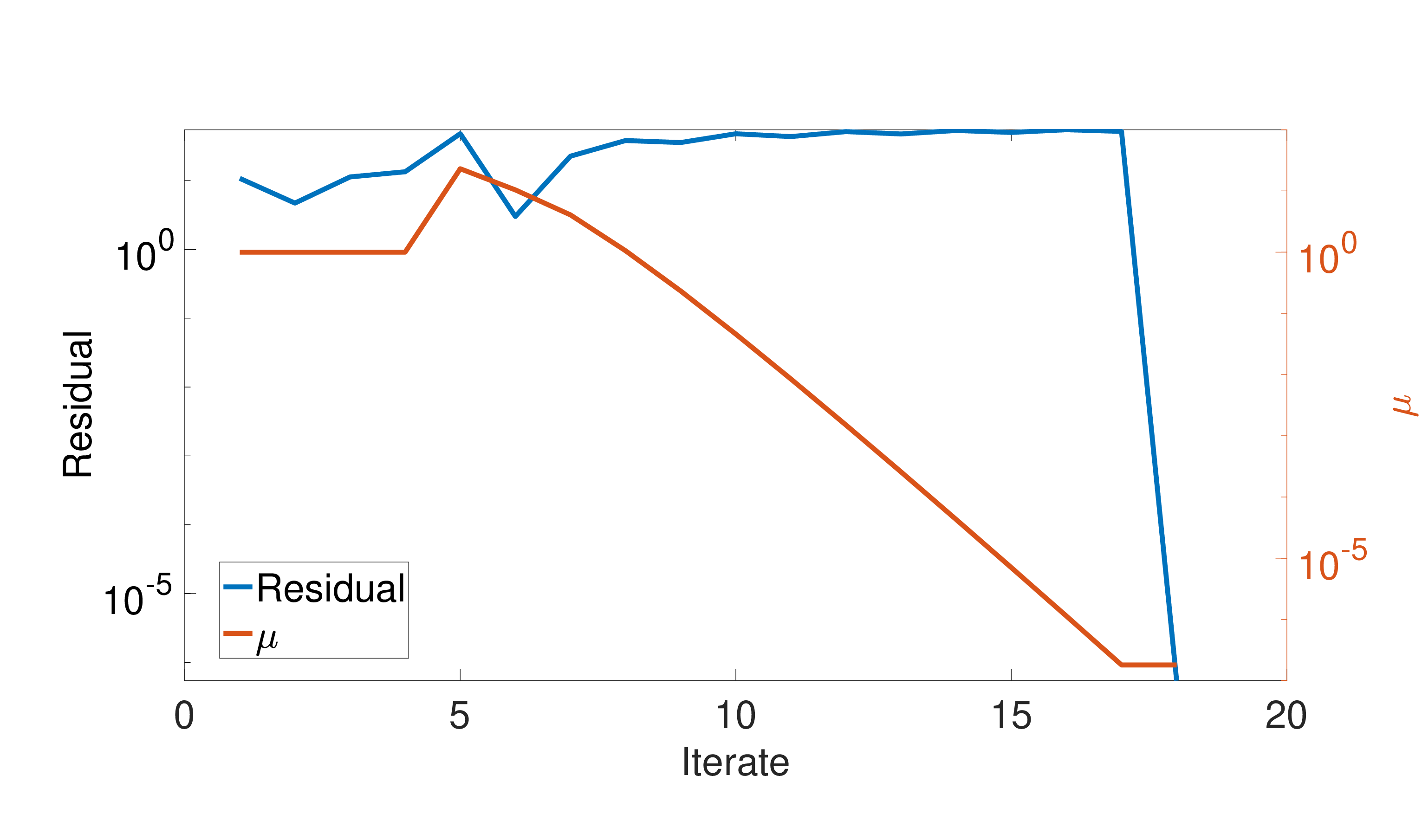}
      \caption{VI residual and barrier parameter $\mu$ at each iteration when solving the sample game. The orange curve denotes the value of $\mu$, and the blue curve denotes the value of residual. It demonstrates that the $\mu$ remains constant during the additional Newton steps while the residual is greatly reduced}
      \label{fig:fig4}
   \end{figure}

As before, we use an sequence $C_i=1+ki$, with $k=0.01$, to ensure the game is non-potential and generate all other parameters randomly. The most important parameter in the game is the sensitivity of the sales price to average sales $C_i$, the player with the largest $C_i$ will have considerably fewer sales than the player with the smallest $C_i$, all else being equal. Fig~\ref{fig:fig3} depicts a representative solution of the market game, we observe that the players with the largest $C_i$ have significantly fewer sales than those with smaller $C_i$. Similar to the previous example, in Fig~\ref{fig:fig4} we observe the same 2-step convergence behaviour with convergence in roughly 20 iterations. The detailed derivation of the KKT system can be found in the Appendix.

\subsection{Example 3: A traffic routing game}
Delays and congestion in traffic networks waste travelers’ time, increase emissions, and reduce urban mobility [3]. Consider the road network modeled as a directed graph $\mathcal{G}=(\mathcal{N},\mathcal{E})$ where $\mathcal{N}\in\{1,...,6\}$ is the set of nodes and $\mathcal{E}\in\{1,...,14\}$ is the set of directed edges that connect the nodes, representing the roads in the traffic network. We model traffic flow as a game where each player wishes to find the fastest route between an origin/destination pair.
   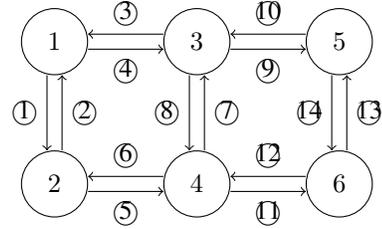
\begin{figure}[thpb]
      \centering
        \tikzset{auto shift/.style={auto=right,->,
            to path={let \p1=(\tikztostart),\p2=(\tikztotarget),
        \n1={atan2(\y2-\y1,\x2-\x1)},\n2={\n1+180}
            in ($(\tikztostart.{\n1})!1mm!270:(\tikztotarget.{\n2})$) -- 
            ($(\tikztotarget.{\n2})!1mm!90:(\tikztostart.{\n1})$) \tikztonodes}}}
        \begin{tikzpicture} 
        \node[state](a) {$1$};
        \node[state,below=of a] (b) {$2$};  
        \node[state,right=of a] (c) {$3$};  
        \node[state,below=of c] (d) {$4$};
        \node[state,right=of c] (e) {$5$};  
        \node[state,below=of e] (f) {$6$}; 
        \draw
        (a) edge[auto shift] node {$\textcircled{1}$} (b)
        (b) edge[auto shift] node {$\textcircled{2}$} (a);  
        \draw
        (a) edge[auto shift] node {$\textcircled{4}$} (c)
        (c) edge[auto shift] node {$\textcircled{3}$} (a);
        \draw
        (b) edge[auto shift] node {$\textcircled{5}$} (d)  
        (d) edge[auto shift] node {$\textcircled{6}$} (b);
        \draw
        (c) edge[auto shift] node {$\textcircled{8}$} (d)
        (d) edge[auto shift] node {$\textcircled{7}$} (c);
        \draw
        (c) edge[auto shift] node {$\textcircled{9}$} (e)
        (e) edge[auto shift] node {$\textcircled{10}$} (c);
        \draw
        (d) edge[auto shift] node {$\textcircled{11}$} (f)
        (f) edge[auto shift] node {$\textcircled{12}$} (d);
        \draw
        (e) edge[auto shift] node {$\textcircled{14}$} (f)
        (f) edge[auto shift] node {$\textcircled{13}$} (e);
        \end{tikzpicture}
      \caption{A road network with 14 edges. Each edge represents a road link with traffic direction and each vertex represents the starting/ending point of the road links}
      \label{fig:network}
   \end{figure}

In this case study, we consider a small traffic system in Fig~\ref{fig:network} with 4 players, denoted by their origin-destination pairs $O_1^5$, $O_1^6$, $O_2^5$, $O_2^6$. Each player represents an origin/destination pair and the $i^{th} player$ has the decision variable $x_i=[x^1_i, x^2_i,...,x^{14}_i]^\intercal$ which is the traffic flow of the $i^{th}$ player on each of the 14 edges in Fig. 5. We define the unit of the decision variable be $\#\,of\,cars/min$ since it represents the traffic flow.

The objective function for the $i^{th}$ player is
\begin{equation}
J_i=\sum_{e}{x_i^e l_e(\mathbf{x}^e)}\tau_i
\end{equation}

where $\tau_i$ is $i^{th}$ player's unit value-of-time ($\$/min$), and $l_e (\mathbf{x}_e)$ is the actual travel time on road indexed by $e$. The actual traffic flow on the road is an aggregative flow of all players in the network so that each player’s routing strategy impacts the congestion outcome for all participants. Let $\mathbf{x}^e=\sum_{e} x_i^e$ be the aggregated traffic flow on road $e$, we can estimate the time a vehicle spends on a road with an affine function of the ratio between the current traffic flow and free-flow traffic flow

\begin{equation}
\begin{aligned}
l_e (\mathbf{x}^e)= t_e(1+4\frac{\mathbf{x}^e}{c_e})\\
c_e=\frac{80}{km} \times (speed\,limit)\times \frac{1}{60}\\
\end{aligned}
\end{equation}

where $t_e$ is the free-flow travel time in minute and $t_e=\frac{road\,length}{speed\,limit}\times 60$, $c_e$ is the capacity of the road indexed by $e$ while we define the length 
of the car be $12.5m/car$. These configurations are the specific conditions of each individual road.

In addition to non-negative constraints for the traffic flow, we enforce a network-flow constraint to ensure all trips are completed. This is accomplished using the affine equality constraint \cite{c14}
\begin{equation}
G_{idc}x_i=\mathbf{s_d^o}(i)
\end{equation}
where $G_{idc}\in \{-1,0,1\}^{|\mathcal{N}|\times |\mathcal{E}|}$ and letting $g_{ie}$ be the element on $i^{th}$ row and $e^{th}$ column of $G_{idc}$

we have 
\begin{equation}
g_{ie} = \begin{cases}
    1 & \text{i is the origin of edge e}\\
    -1 & \text{i is the destination of edge e}\\
    0 & \text{otherwise}\\
\end{cases}
\end{equation}

$\mathbf{s_d^o}(i)$ is the source-sink vector for $i^{th}$ player. We use ${s_d^o}^i$ to denote the total traffic flow for the $i^{th}$ player, which is the sum of all the traffic flows among all edges belonging to the $i^{th}$ player, with the origin being $o$ and destination being $d$. If we define $j$ be index in the ${s_d^o}(i)$ vector then

\begin{equation}
\mathbf{s_d^o} = \begin{cases}
    s_d^o & j =o\\
    -s_d^o & j=d\\
    0 & otherwise.
\end{cases}
\end{equation}

There are no coupling constraints, the local feasible set is simply
\begin{equation}
\begin{aligned}
\{x_i\in \mathbb{R}^{14}| x^i \geq \mathbf{0}, \,G_{idc}x_i=\mathbf{s_d^o}(i),\, \forall i\in[4]\}
\end{aligned}
\end{equation}

The detailed KKT system derivation for this problem is in the Section C of the Appendix.

 In this problem the parameters ${\tau_1,...,\tau_4}$ determine if the VI is symmetric, as in the previous example we generate the parameters using the equation using the sequence $\tau_i=1+ki$, with $k=0.1$ to generate a non-symmetric problem. The road lengths (in km) are $l=[2,\, 2.5,\, 1,\, 5,\, 3,\, 5,\, 1.5,\, 6,\, 3,\, 5,\, 4,\, 2,\, 1,\, 2]$ for road indexed from 1 to 14, with speed limits (in km/h) being $s=[50,\, 50,\, 30,\, 50,\, 50,\, 80,\, 30,\, 80,\, 50,\, 50,\, 80,\, 50,\, 30,\, 50]$.

\begin{figure}[thpb]
      \centering
      \includegraphics[scale=0.2]{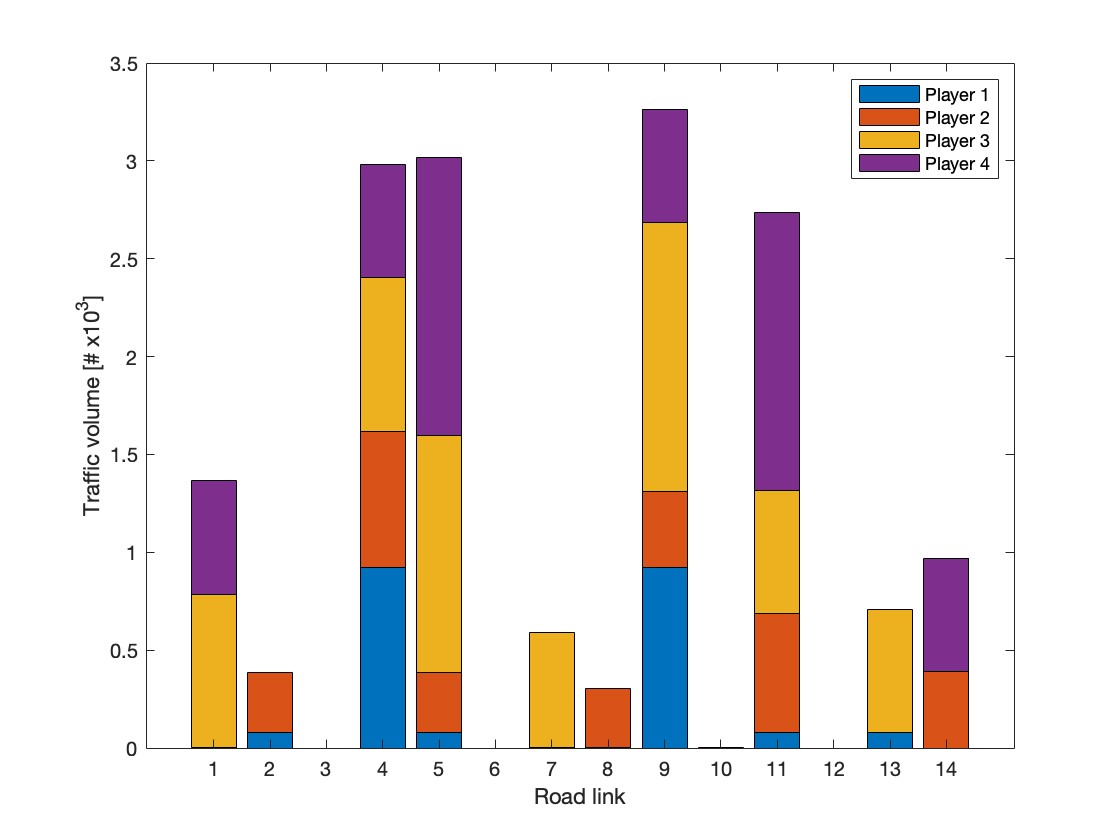}
      \caption{Demonstration of the solution to the randomly generated vehicle routing game. The distribution of traffic flow on the road network is shown on the stacked bar plot.}
      \label{fig:fig5}
   \end{figure}

\begin{figure}[thpb]
      \centering
      \includegraphics[scale=0.15]{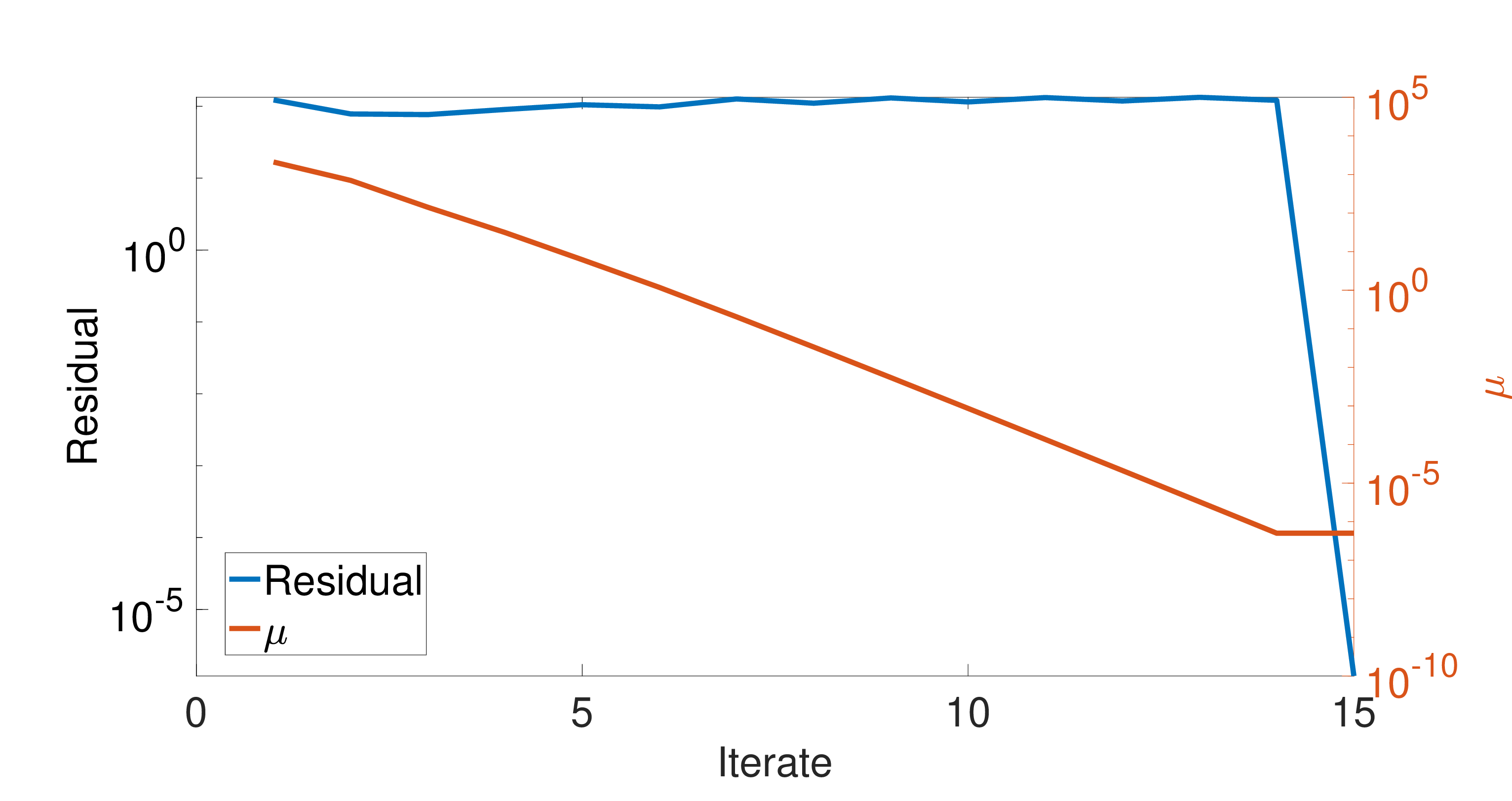}
      \caption{VI residual and barrier parameter $\mu$ at each iteration when solving the sample problem. The orange curve denotes the value of $\mu$, and the blue curve denotes the value of residual. It demonstrates that the $\mu$ remains constant during the additional Newton steps.}
      \label{fig:fig7}
   \end{figure}

Fig~\ref{fig:fig5} demonstrates the solution to the problem. We find that a few road links have been avoided by all traffic, which is reasonable since they are in reverse directions to the destination and will certainly increase the time delay if chosen. We also determine that certain groups of users primarily use a few roads. Similar to the previous examples, in Fig~\ref{fig:fig7}, we found the residual appears to have a 2-stage behavior and converges quickly with low residual.

\section{Numerical Benchmarking}
We perform a simulation study to investigate the impact of problem size and monotonicity constant on solver performance. We benchmark our IPM against the Forward-backward Splitting (FBS) method proposed in \cite{c17} alongside the EV charging benchmark problem. All the tests are run on a MacBook Air laptop with Apple M2 chip.

\subsection{Performance vs. Problem Size}
In this section, we investigate how the performance of our method scales with problem size. We generate a collection of benchmark problems based on market game and EV charging game as described in the previous sections while varying the number of players in the game from 5 to 50 in increments of five (keeping the number of markets fixed at $m = 10$ and the number of hours fixed at $n = 24$).

For each problem size, we use both methods to solve the same problem ten times and record the average wall clock time The convergence condition for FBS methods is set to be $\|x^k - x^*\|/\|x^*\| \leq10^{-6}$ as in \cite{c9} where $x^*$ denotes a high-precision solution of \eqref{eq:VI} and $x^k$ denotes the current iteration. The tolerances of the IPM were set as $\varepsilon = \tau = 10^{-6}$ to ensure that the precision of the IPM and FBS solutions were comparable.

 \begin{figure}[thpb]
      \centering
      \includegraphics[scale=0.15]{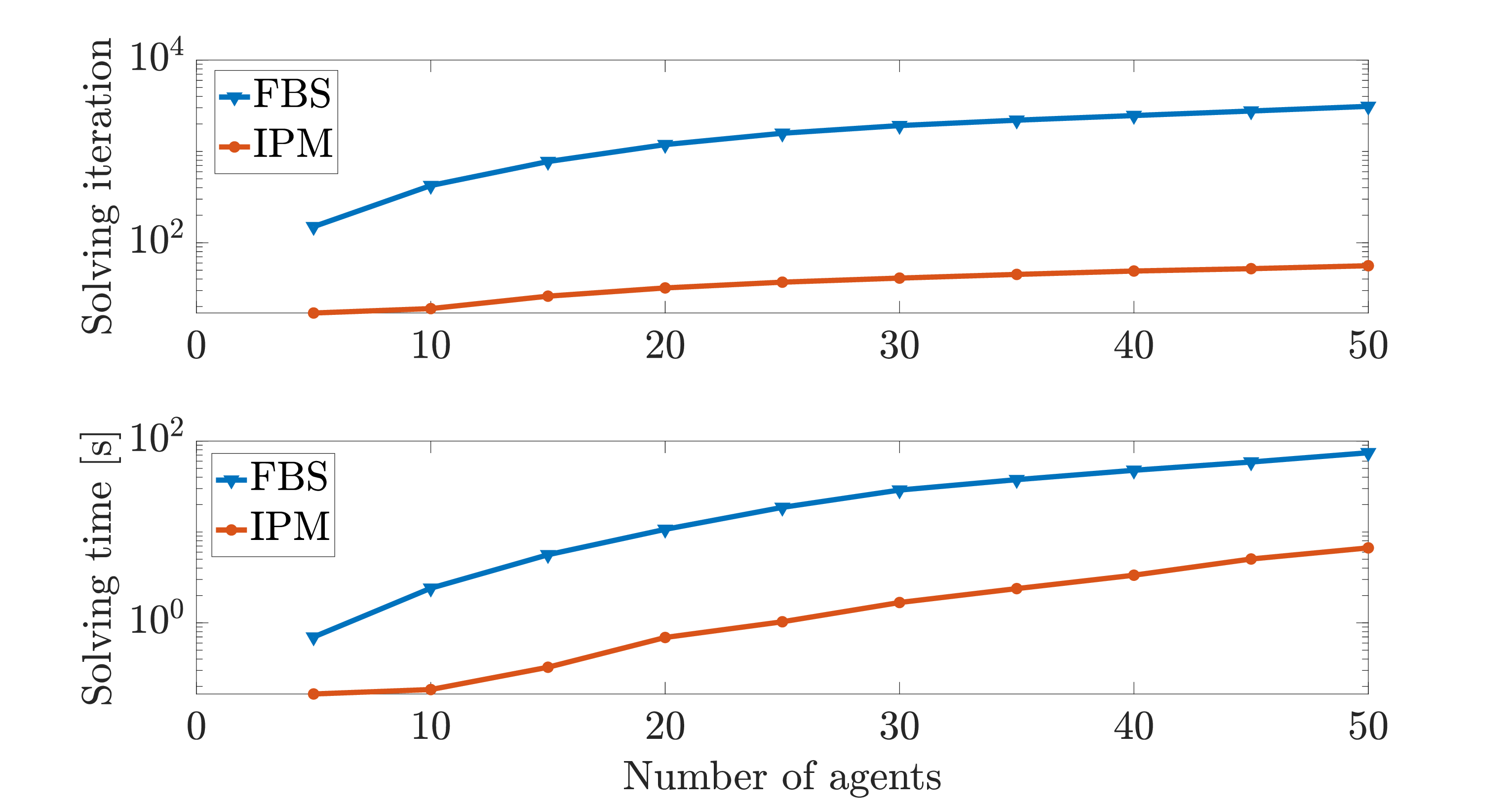}
      \caption{Iteration count and solving time vs. different problem sizes when solving the market game with the FBS and IPM methods.}
      \label{fig:fig8}
   \end{figure}
   
  \begin{figure}[thpb]
      \centering
      \includegraphics[scale=0.15]{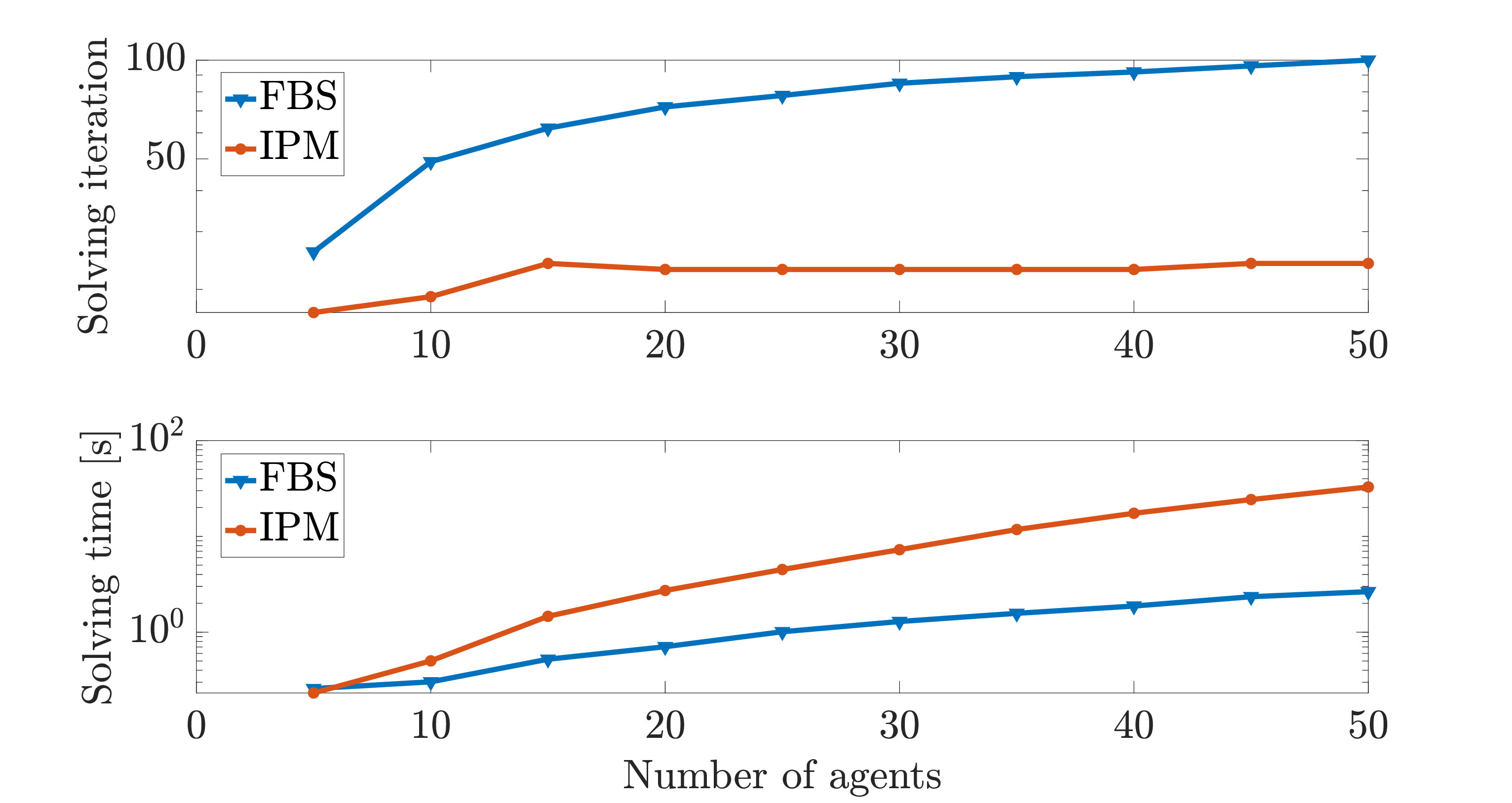}
      \caption{Iteration count and solving time vs. different problem sizes when solving the EV charging game with the FBS and IPM methods.}
      \label{fig:fig9}
   \end{figure}
   
Fig~\ref{fig:fig8} and Fig~\ref{fig:fig9} compare the iteration count and execution time for the two methods. The IPM converges in far fewer iterations and the number of iterations grows slowly (essentially plateauing) for both problems as expected for a second-order method. In terms of time, the IPM outperforms the FBS for the market game but will eventually be slower for very large-size problems. The FBS method requires very few iterations to converge for the EV charging game and beats the IPM. This suggests that second-order methods are best suited for ``harder'' (i.e., poorly conditioned problems where first-order methods will converge slowly) of small-medium scales (since second-order methods are difficult to parallelize/distribute).

%  Overall, the FBS method is more sensitive to the problem setup, performance of the IPM is stable and can solve the problems in a reasonably short time among the market game and EV charging game. However, the FBS method takes a much longer time to solve the market game potentially due to its more complex constraints.

% Fig. 8 shows that 

% Furthermore, Fig. 9 shows that for a simpler problem (EV charging game), the FBS method outperforms IPM in both ways but the difference is small in terms of the solving time. We find that the IPM method's performance is more stable and we determine that the IPM method is best for small to medium size games.

% In fact, when solving QG with simpler constraints, such as the EV charging game, the FBS method is faster in most cases while IPM is still highly efficient. We determine that the FBS method's step size is sensitive to the constraints. For more complex problems, the FBS method will not only be slower but require much more effort in reformulating the solver to fit the new problem. For instance, it can be considerably onerous to apply the FBS method to the traffic routing problem while its performance is likely not comparable to IPM.

\subsection{Solver performance with problem monotonicity}
The monotonicity constant $\lambda_{\min} (0.5(W + W^T))$ of the game tends to have a strong impact on the convergence rate of FBS and other first-order algorithms. In this section, we investigate the impact of the monotonicity constant on the performance of the FBS and IPM algorithms. 

We randomly generate EV charging and Market games of a fixed-size as described in Section~\ref{ss:case-studies} and adjust the parameters $\{C_1,C_2,...,C_N\}$ (which govern the montonicity properties of both games) by setting $C_i = 1 + ki$ and varying $k$ from 0.1 to 3. This generated a set of problems with monotonicity constants ranging from 0.2 to 0.7. We solved each problem instance 10 times and averaged the execution time results.

% We set the game to ten agents competing in ten markets for the market problem. We increase the k value from 0.1 to 3 to generate random problems with 

% The monotonicity constant of a QG is defined by $Definition$ $1$, and we manually define $\{C_1,C_2,...,C_N\}$ to be an sequence such that $C_i=1+ki$. We can tune $k$ to adjust the problem's monotonicity constant to qualitatively show its impacts on the solver's performance.

\begin{figure}[thpb]
      \centering
      \includegraphics[scale=0.15]{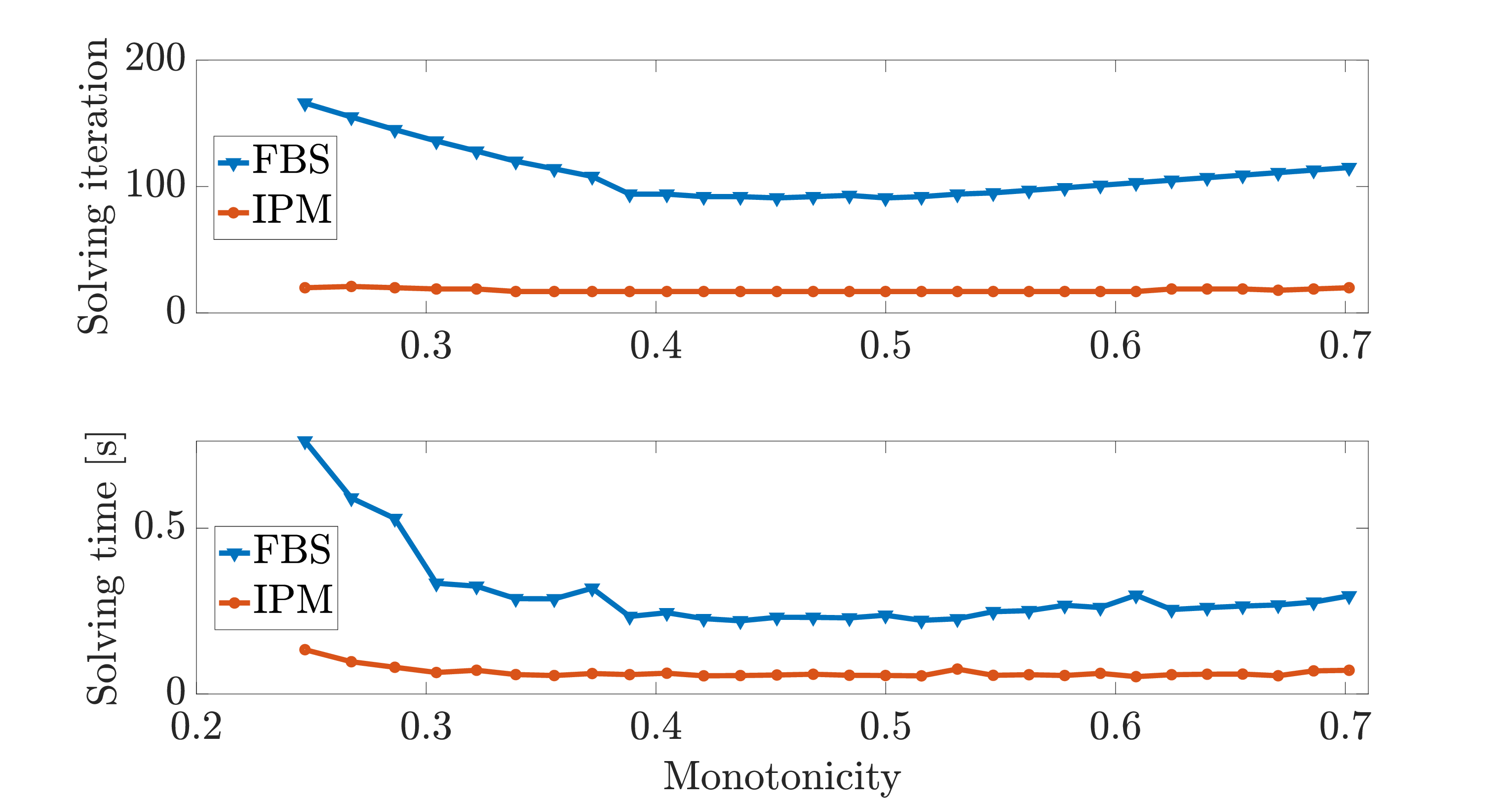}
      \caption{Iteration count and solving time vs. different problem monotonicity constants when solving the market game with the FBS and IPM methods.}
      \label{fig:fig10}
   \end{figure}
   
    \begin{figure}[thpb]
      \centering
      \includegraphics[scale=0.15]{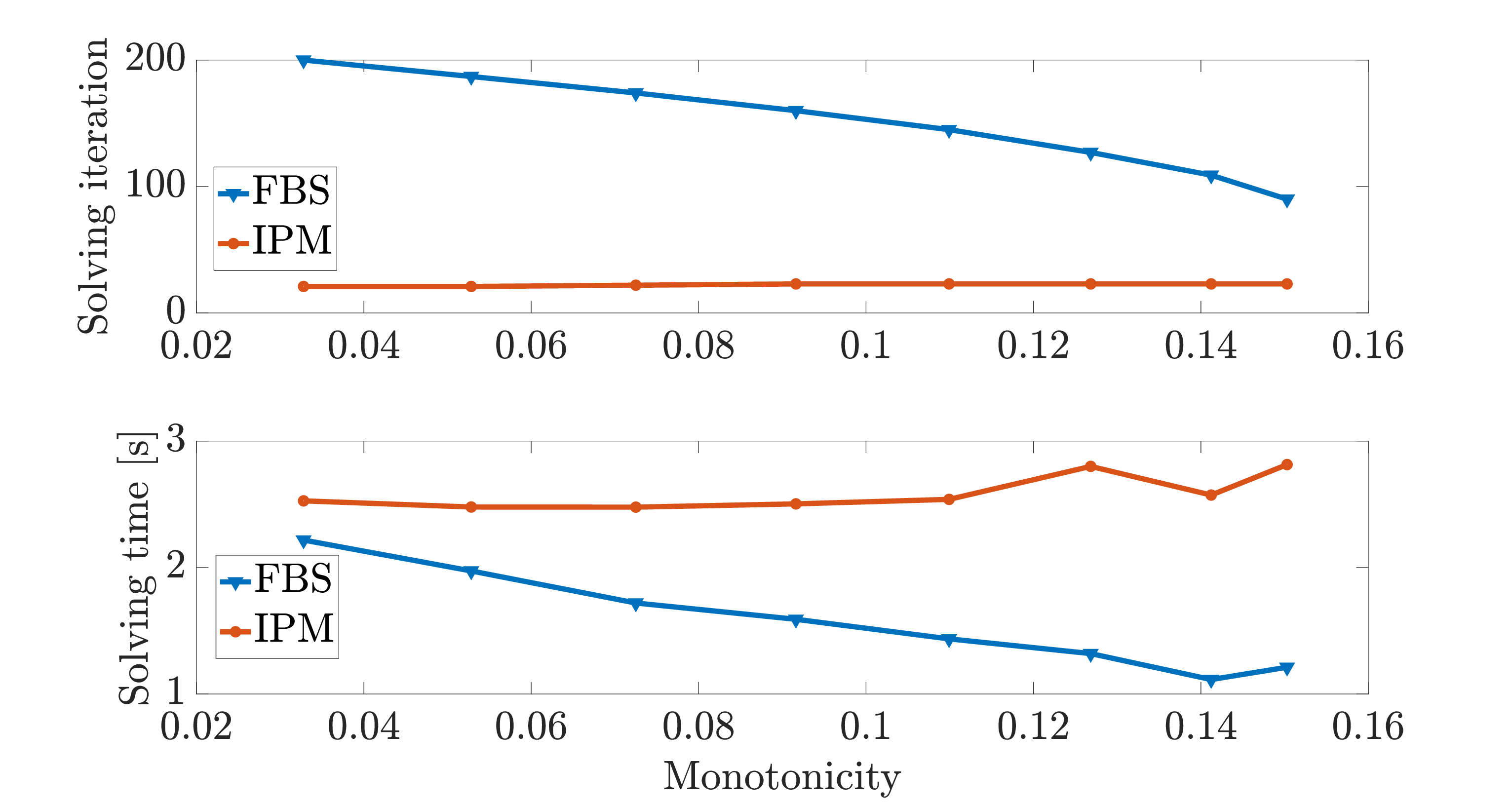}
      \caption{Iteration count and solving time vs. different monotonicity constants when solving the EV charging game with the FBS and IPM methods.}
      \label{fig:fig11}
   \end{figure}

The results are shown in Figures~\ref{fig:fig10} and \ref{fig:fig11}. The FBS method is faster when the monotonicity constant increases (up to some limit) while it does not impact the performance of the IPM.

Overall, we found that the IPM is a stable method for small to medium size problems that is less sensitive to variations in the problem properties than the FBS method. The log-domain IPM is a promising algorithm for constructing a general purpose solver for CQGs.

% It appears that the FBS method becomes faster at problems with a higher monotonicity constant, as shown in Fig 9 and Fig 10. We determine that the monotonicity constant of the problem does not influence the performance of IPM significantly. We determine that IPM is not only faster but more stable when the problem's monotonicity constant changes.
 
\section{CONCLUSIONS}
This paper presents an algorithms for finding GNEs of CQGs based on a newly proposed log-domain interior point method. This IPM is easy to adapt to different problems, relatively insensitive to problem properties (e.g., the monotonicity constant) and scales well for small to medium size problems, making it attractive as a general-purpose algorithm for finding v-GNEs of convex-quadratic games. Future work includes developing a higher performance open-source implementation and validating its performance on a wider array of test problems.

\section*{ACKNOWLEDGMENT}
The authors would like to thank Dr. Frank Permenter for providing technical advice at the early stage of this study, Dr. Giuseppe Belgioioso for providing test scripts.

%%%%%%%%%%%%%%%%%%%%%%%%%%%%%%%%%%%%%%%%%%%%%%%%%%%%%%%%%%%%%%%%%%%%%%%%%%%%%%%%

%                                   % the textheight of the last page by a suitable amount.
%                                   % This command does not take effect until the next page
%                                   % so it should come on the page before the last. Make
%                                   % sure that you do not shorten the textheight too much.

% \clearpage
\appendix
\subsection{Electric vehicle charging game}
The target is to derive the concatenated KKT system into the format of (4). We first derive both the individual and coupling constraints into an unified affine format such as $A\mathbf{x}+b\geq \mathbf{0}$.

We express the individual constraint for the $i^{th}$ player as $A_i x_i+b_i \geq \mathbf{0}$ where
\begin{equation}
A_i=
\begin{bmatrix}
1 &  & \\
& \ddots &\\
&  &  1 \\
-1 &  & \\
& \ddots & \\
&  & -1\\
1&    \hdots &     1\\
\end{bmatrix}
\end{equation}\\
and
\begin{equation}
b_i=
\begin{bmatrix}
0\\
\vdots\\
0 \\
0\\
\vdots\\
0 \\
-l_i\\
\end{bmatrix}
\end{equation}\\
Matrix $A_i$ implements $(13)$, and the last row of matrix $A_i$ is to implement that the charging level is satisfactory at the end of the time period.

We express the coupling constraints as $A_{couple}\mathbf{x}+NK(t)\geq \mathbf{0}$ where

\begin{equation}
A_{couple}= {\begin{bmatrix}
{\begin{bmatrix}
-1 & &\\
& \ddots &\\
& & -1\\
\end{bmatrix}}
\cdots
{\begin{bmatrix}
-1 & &\\
& \ddots &\\
& & -1\\
\end{bmatrix}}\\
\end{bmatrix}}
\end{equation}
and $K$ is a vector of size n while $t\in[n]$.
We concatenate the constraints, and derive the overall feasible set $\chi=\{\mathbf{x} \in \mathbb{R}^{nN}:A\mathbf{x}+b\geq \mathbf{0}\}$, where
\begin{equation}
A= 
{\begin{bmatrix}
A_1 & &\\
& \ddots &\\
& & A_N\\
A_{couple} & &\\
\end{bmatrix}}
\end{equation}
and 
\begin{equation}
b= 
{\begin{bmatrix}
b_1\\
\vdots\\
b_N\\
NK_t\\
\end{bmatrix}}
\end{equation}

We then derive the concatenated pseudo-gradient of the objective functions. We first derive the pseudo-gradient of the objective function of $i^{th}$ player.
\begin{equation}
\nabla_{x_i} J_i(x_i, \Bar{\textbf{x}})=Q_i^\intercal x_i + p_i + C_i \nabla_{x_i}(\Bar{\textbf{x}}^\intercal) x_i + (C_i \Bar{\textbf{x}}) + c)^\intercal
\end{equation}

where\\
\begin{equation}
 Q_i=
 {\begin{bmatrix}
q(1) & &\\
&  \ddots &\\
&    &   q(n)\\
\end{bmatrix}}
\end{equation}
and\\
\begin{equation}
 \nabla_{x_i}(\Bar{\textbf{x}}^\intercal)=
 {\begin{bmatrix}
\frac{1}{N} & &\\
&  \ddots &\\
&    &   \frac{1}{N}\\
\end{bmatrix}}
\end{equation}
also
\begin{equation}
\Bar{\textbf{x}}=
 {\begin{bmatrix}
\frac{1}{N} & & &\frac{1}{N} & \hdots  & &\frac{1}{N} & &\\
& \ddots & & & \ddots &  & & \ddots &\\
& & \frac{1}{N} & & & \frac{1}{N} & \hdots & &\frac{1}{N}\\
\end{bmatrix}}
\textbf{x}
\end{equation}
We concatenate the pseudo-gradient of the game according to Theorem 3.3 of \cite{c9} into the form of \eqref{eq:W-def} where
\begin{equation}
W=
{\begin{bmatrix}
W_1 &  &  &\\
&  W_2 &  &\\
&   &  \ddots  &\\
&   &  &  W_N\\
\end{bmatrix}} + 
{\begin{bmatrix}
C_1\,avg\\
C_2\,avg\\
\vdots\\
C_N\,avg\\
\end{bmatrix}}
\end{equation}
of which
\begin{equation}
W_i=
Q_i+C_i\nabla_{x_i}(\Bar{\textbf{x}}^\intercal) 
\end{equation}
and
\begin{equation}
avg=
 {\begin{bmatrix}
\frac{1}{N} & & &\frac{1}{N} & \hdots  & &\frac{1}{N} & &\\
& \ddots & & & \ddots &  & & \ddots &\\
& & \frac{1}{N} & & & \frac{1}{N} & \hdots & &\frac{1}{N}\\
\end{bmatrix}}
\end{equation}
and 
\begin{equation}
f= {\begin{bmatrix}
p_1 +c\\
\vdots\\
p_N+c\\
\end{bmatrix}}
\end{equation}
while $c=(d(1),...,d(n))$ as we assume all EVs are faced with the same non-EV background power demand.
\subsection{Market game}
We first derive the constraints of the market game into an affine format. For $i^{th}$ player, we express its individual constraint as:
\begin{equation}
A_i x_i+b_i \geq 0
\end{equation}
where 
\begin{equation}
A_i= {
\begin{bmatrix}
1 & & & & &\\
& \ddots & & & &\\
& & 1 & & &\\
& & & 1 & &\\
& & & & \ddots &\\
& & & & & 1\\
-1 & & & & &\\
& \ddots & & & &\\
& & -1 & & &\\
1 & \cdots & 1 &-1 & \cdots & -1\\
\end{bmatrix}}
\end{equation}
and
\begin{equation}
b_i= {
\begin{bmatrix}
0\\
\vdots\\
0\\
0\\
\vdots\\
0\\
u_{i,1}\\
\vdots\\
u_{i,m}\\
0\\
\end{bmatrix}.}
\end{equation}

Similarly, we express the coupling constraints into the form of $A_{couple}\mathbf{x}+b_{couple} \geq 0$
where 
\begin{equation}
A_{couple} = {\begin{bmatrix}
I_m0_m
\cdots
I_m0_m\\
-I_m0_m
\cdots
-I_m0_m
\end{bmatrix}}
\end{equation}
where
\begin{equation}
I_m0_m = 
{\begin{bmatrix}
1 & & & 0 & &\\
& \ddots & & & \ddots &\\
& & 1 & &  &0\\
\end{bmatrix}}
\end{equation}
and
\begin{equation}
b_{couple}= {
\begin{bmatrix}
0\\
\vdots\\
0\\
0\\
\vdots\\
0\\
\end{bmatrix}}
\end{equation}

Then the overall constraint for the game in affine format is $A \mathbf{x} +b \geq 0$ where,
\begin{equation}
A=
{\begin{bmatrix}
A_1  &  &  &\\
&    A_2  &  &\\
&   &  \ddots  &\\
&   &     &    A_N\\
A_{couple} &  &  & \\
\end{bmatrix}}\\
\end{equation}
and
\begin{equation}
b=
{\begin{bmatrix}
b_1\\
b_2\\
\vdots\\
b_N\\
b_{couple}\\
\end{bmatrix}}\\
\end{equation}
The overall affine feasible set is $\mathbf{\chi}=\{\mathbf{x}\in\mathbb{R}^{nN}| A\mathbf{x}+b\geq 0\}$.

Then we express the pseudo-gradient for the objective function of $i^{th}$ player
\begin{equation}
\nabla_{x_i} J_i(x_i, \bar{\mathbf{x}}) = 2Q_i x_i + p_i^\intercal + C_i(\nabla_{x_i}(\Delta \bar{\mathbf{x}})^\intercal) x_i + C_i\Delta \bar{\mathbf{x}}
\end{equation}
where
\begin{equation}
\nabla_{x_i}(\Delta \bar{\mathbf{x}})=
\begin{bmatrix}
0 & & & & &\\
& \ddots & & & &\\
& & 0 & & &\\
& & & 1/N & &\\
& & & & \ddots &\\
& & & & & 1/N\\
\end{bmatrix}.
\end{equation}
Also, we can express
\begin{equation}
\Delta \bar{\mathbf{x}}=G\mathbf{x}
\end{equation}
where $G=[\tilde{G} \cdots \tilde{G}]$
\begin{equation}
\tilde{G}=
\begin{bmatrix}
0 & & & & &\\
& \ddots & & & &\\
& & 0 & & &\\
& & & 1/N & &\\
& & & & \ddots &\\
& & & & & 1/N\\
\end{bmatrix}
\end{equation}
We concatenate the pseudo-gradient of the game according in the format of $(6)$
in which
\begin{equation}
W=
\begin{bmatrix}
W_1 & & \\
& \ddots &\\
& & W_N\\
\end{bmatrix}
+
\begin{bmatrix}
C_1G\\
\vdots \\
C_NG\\
\end{bmatrix}
\end{equation}
and $W_i=2Q_i^\intercal +C_i(\nabla_{x_i}(\Delta \bar{\mathbf{x}})^\intercal)$.

\subsection{Traffic routing game}
Let $\mathbf{x}=(x_1,x_2,x_3,x_4)$, the overall feasible set is then 
\begin{equation}
\begin{aligned}
\mathbf{\Omega}=\{\mathbf{x}\in \mathbb{R}^{64}| A\mathbf{x}+b \geq \mathbf{0}, \,G\mathbf{x}+h=\mathbf{0},\, \forall i\in[4]\}
\end{aligned}
\end{equation}
in which, for this specific problem setup,
\begin{equation}
\begin{aligned}
A=I_{14\times4},
b=\mathbf{0}
\end{aligned}
\end{equation}
and
\begin{equation}
G=\begin{bmatrix}
G_{idc} & & & \\
& G_{idc} & &\\
& & G_{idc} &\\
& & & G_{idc}\\
\end{bmatrix},
h=\begin{bmatrix}
-\mathbf{s_d^o}(1)\\
-\mathbf{s_d^o}(2)\\
-\mathbf{s_d^o}(3)\\
-\mathbf{s_d^o}(4)\\
\end{bmatrix}
\end{equation}
The pseudo-gradient for player 1, $O_1^5$ is denoted below
\begin{equation}
 \nabla_{x_1}(J_1)=
 {\begin{bmatrix}
2W_0\tau_1 & W_0\tau_1& W_0\tau_1& W_0\tau_1\\
\end{bmatrix}}
\textbf{x}+\tau_1\mathbf{t}
\end{equation}
where
\begin{equation}
 W_0=
 {\begin{bmatrix}
\frac{4}{c_1} \times 2t_1 & &\\
& \ddots & \\
& & \frac{4}{c_{14}}\times 2t_{14}\\
\end{bmatrix}},
\mathbf{t}={\begin{bmatrix}
t_1\\
\vdots\\
t_{14}\\
\end{bmatrix}}
\end{equation}
We concatenate the pseudo-gradient of the game according to Theorem 3.3 of [5] in the format of $(6)$, in which
\begin{equation}
W=\begin{bmatrix}
2W_0\tau_1 & W_0\tau_1& W_0\tau_1& W_0\tau_1\\
W_0\tau_2 & 2W_0\tau_2& W_0\tau_2& W_0\tau_2\\
W_0\tau_3 & W_0\tau_3& 2W_0\tau_3& W_0\tau_3\\
W_0\tau_4 & W_0\tau_4& W_0\tau_4& 2W_0\tau_4\\
\end{bmatrix},
f=\begin{bmatrix}
\tau_1\mathbf{t}\\
\tau_2\mathbf{t}\\
\tau_3\mathbf{t}\\
\tau_4\mathbf{t}\\
\end{bmatrix}
\end{equation}
We then convert the game into its VI and could then solve the VI with the IPM solver.

% 

%                                   % the textheight of the last page by a suitable amount.
%                                   % This command does not take effect until the next page
%                                   % so it should come on the page before the last. Make
%                                   % sure that you do not shorten the textheight too much.

\end{document}